\documentclass[12pt,leqno]{article}

\usepackage{amssymb}
\usepackage{amsmath}
\usepackage{graphicx,psfrag,color}
\usepackage{url}
\usepackage{pstool}

\newcommand{\comment}[1]{}

\newcommand{\dm}[1]{\mathrm{dm(#1)}}
\newcommand{\rg}[1]{\mathrm{rg(#1)}}

\newcommand{\conv}{{}^{\smallsmile}}
\newcommand{\inv}{{}^{-1}}
\newcommand{\rp}{\mathrel{\vert}}

\newcommand{\id}{\mbox{\sf Id}}

\newcommand{\Z}{\mathcal{S}}
\renewcommand{\P}{F}
\newcommand{\HH}{\mathcal{H}}
\newcommand{\rep}{\mbox{\sf rep}}

\newcommand{\de}{\mbox{\ :=\ }}
\newcommand{\deiff}{\mbox{$\quad : \Leftrightarrow\quad$}}
\newcommand{\Aa}{\mbox{$\mathfrak A$}}
\newcommand{\qed}{\hfill\mbox{\bf QED}\bigskip}

\newtheorem{thm}{Theorem}

\newtheorem{rem}{Remark}

\begin{document}

\title{Term algebras of elementarily equivalent atom structures}
\author{Andr\'eka, H. and N\'emeti, I.}%
\date{}
\maketitle

\centerline{We dedicate this paper to Bjarni J\'onsson.}

\begin{abstract} We exhibit two relation algebra atom structures such that they are
elementarily equivalent but their term algebras are not. This
answers Problem 14.19 in the book Hirsch, R.\ and Hodkinson, I.,
``Relation Algebras by Games", North--Holland, 2002.
\end{abstract}

\section{Introduction}
Atom structures for Boolean algebras with operators, and in
particular for relation algebras, were introduced in
J\'onsson-Tarski~\cite{JT}. These structures proved to be a central
tool in algebraic logic, see, e.g., \cite{CAIII},
\cite[section~19]{G17}, \cite[section~2.7]{HMT71},
\cite[section~2.7]{HH02}, \cite{MTAMS, Mbook, SDis}. There is a kind
of duality between atom structures of algebras and complex algebras
of relational structures, a large part of this duality is elaborated
in \cite{JT}. Atom structures are useful in constructing relation
algebras, because atom structures are simpler and hence easier to
work with. Therefore, it is useful to see what properties of complex
algebras can be ensured by constructing appropriate atom structures.
Representability of the complex algebra is not such a property,
because there are two elementarily equivalent relation algebra atom
structures such that the complex algebra of one is representable
while that of the other is not, a result of Hirsch and
Hodkinson~\cite[Corollary~14.14]{HH02}.

The term algebra of an atom structure is the smallest subalgebra of
its complex algebra that has the same atom structure: it is the
subalgebra of the complex algebra generated by the singletons of the
atoms. They are more tightly connected to their atom structures,
e.g., representability of the term algebra of a relation algebra
atom structure can be ensured by its first-order logic theory. So,
if two atom structures are elementarily equivalent, then their term
algebras are either both representable, or neither of them is
representable (a result of Venema~\cite{Ven}, see
\cite[Theorem~2.84]{HH02}). If two atom structures are ``very
close", i.e., if they are $L_{\infty,\omega}$-indistinguishable,
then their term algebras are also this very close, i.e., they are
$L_{\infty,\omega}$-indistinguishable (\cite[Exercise 14.6]{HH02}).
Problem 14.19 in \cite{HH02} asks if this last result holds with
$L_{\omega\omega}$ in place of $L_{\infty,\omega}$: If $S, S'$ are
elementarily equivalent relation algebra atom structures, must the
term algebras of $S$ and $S'$ also be elementarily equivalent?

In this paper we give a negative answer to this question. We
construct a relation algebra atom structure $S$, such that the term
algebra of $S$ and that of an ultraproduct of $S$ are not
elementarily equivalent. Moreover, $S$ is an atom structure of a
completely representable simple relation set algebra.

\section{The construction}

We begin by recalling terminology from \cite{HH02}. A (relational
algebra type) \emph{atom structure} is a structure $\langle S, P, C,
I\rangle$, where $P, C$ and $I$ are ternary, binary, and unary
relations on $S$, respectively.  The \emph{complex algebra} of an
atom structure $\langle S,P,C,I\rangle$ is the algebra $\langle A,+,
- , ;,\conv,1'\rangle$ where $A$ is the collection of all subsets of
$S$, $+$ and $-$ are the operations of forming union and complement
(with respect to $S$) respectively, and the operations $;,\conv,1'$
are determined by $P, C, I$ as follows. Let $X,Y\subseteq S$. Then
$X;Y=\{ u : P(x,y,u) \text{ for some } x\in X, y\in Y\}$, $X\conv=\{
u : C(x,u) \text{ for some }x\in X\}$, and $1'=I$. The atom
structure is called \emph{completely representable} if its complex
algebra is completely representable, that is to say if, up to an
isomorphism, $A$ is a set of binary relations such that the biggest
element is the union (as opposed to the supremum only) of the atomic
relations, and the operations $+, -, ;, \conv, 1'$ are the following
standard operations on binary relations: union, complement (with
respect to a largest element of $A$), relation composition,
converse, and the identity relation. Finally, the \emph{term
algebra} of an atom structure is the subalgebra of its complex
algebra generated by the singletons. In the paper, $\omega$ denotes
the set of non-negative integers.

\begin{thm}\label{main-thm}
There are completely representable relational algebra atom
structures which are elementarily equivalent but their term algebras
are not elementarily equivalent.
\end{thm}
\noindent{\bf Proof.}
First we define a relation algebra type atom structure $\Z=\langle
S,P,C,I\rangle$. The universe $S$ of the atom structure is
\begin{description}
\item{}
$\{\id_{n,i} : n\in\omega, i\le n\}\cup$
\item{}
$\{ r_{n,k} : n\in\omega, 1\le k\le n\}\cup \{ r_{n,k}^- :
n\in\omega, 1\le k\le n\}\cup$
\item{}
$ \{ w_{n,i,m,j} : n,m\in\omega, i\le n, j\le m\} ,$
\end{description}
and $I=\{ \id_{n,i} : n\in\omega, i\le n\}$. The binary relation $C$
is
\begin{description}
\item{}
$\{ (x,x) : x\in I\}\cup$
\item{}
$\{ (r_{n,k},r^-_{n,k}) : n\in\omega, 1\le k\le n\}\cup
\{ (r^-_{n,k},r_{n,k}) : n\in\omega, 1\le k\le n\}\cup$
\item{}
$\{ (w_{n,i,m,j}, w_{m,j,n,i}) : n,m\in\omega, i\le n, j\le m\}$.
\end{description}
To define $P$, we first define two unary operations on $S$, the
domain $\dm{a}$ and the range $\rg{a}$ of $a$:
\begin{description}
\item{} $\dm{x}=\rg{x}=x$ for each $x\in I$,
\item{} $\dm{r_{n,k}}=\rg{r^-_{n,k}}=\id_{n,0}$,\quad $\rg{r_{n,k}}=\dm{r^-_{n,k}}=\id_{n,k}$,\quad and
\item{} $\dm{w_{n,i,m,j}}=\rg{w_{m,j,n,i}}=\id_{n,i}$.
\end{description}
Now, the ternary relation $P$ on $S$ is:
\begin{description}
\item{} $\{ (\dm{a},a,a) : a\in S\}\cup \{(a,\rg{a}, a) : a\in
S\}\cup\{ (a,b,\dm{a}) : (a,b)\in C\}\cup$
\item{} $\{ (a,b,c) : a,b,c\in S\setminus I, \dm{a}=\dm{c}, \rg{a}=\dm{b}, \rg{b}=\rg{c}\}$.
\end{description}
With this, the atom structure $\Z$ has been defined.

\begin{figure}[h!]
\begin{center}
\psfragfig[width=\textwidth]{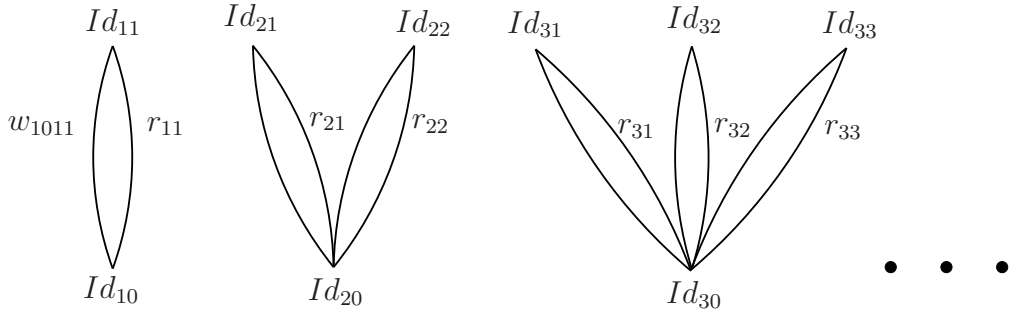} {
\psfrag{a}[b][b][1]{$Id_{11}$} \psfrag{b}[t][t][1]{$Id_{10}$}
\psfrag{c}[r][r][1]{$w_{1011}$} \psfrag{d}[l][l][1]{$r_{11}$}
\psfrag{e}[b][b][1]{$Id_{21}$} \psfrag{f}[b][b][1]{$Id_{22}$}
\psfrag{g}[b][b][1]{$Id_{31}$} \psfrag{h}[b][b][1]{$Id_{32}$}
\psfrag{i}[b][b][1]{$Id_{33}$} \psfrag{j}[l][l][1]{$r_{21}$}
\psfrag{k}[l][l][1]{$r_{22}$} \psfrag{l}[l][l][1]{$r_{31}$}
\psfrag{m}[l][l][1]{$r_{32}$} \psfrag{n}[l][l][1]{$r_{33}$}
\psfrag{o}[t][t][1]{$Id_{20}$} \psfrag{p}[t][t][1]{$Id_{30}$} }
\end{center}
\caption{\label{glove-fig} The split atoms in the atom structure
$\Z$.}
\end{figure}

Here is a short intuitive description of $\Z$, see
Figure~\ref{glove-fig}.
Keeping the complex algebra in mind, we call the elements of $S$
atoms, the elements of $I$ identity atoms, and the rest diversity
atoms, and we say that $a$ goes from $\dm{a}$ to $\rg{a}$. Now,
between any two identity atoms in $\Z$ only one or two diversity
atoms go. The pairs of identity atoms that have two diversity atoms
between them (in both directions) are the pairs $\id_{n,0}$ and
$\id_{n,k}$ with $1\le k\le n$, we call these \emph{splitable pairs}
and the  atoms going between them are called \emph{split} atoms:
from $\id_{n,0}$ to $\id_{n,k}$ the two diversity atoms $r_{n,k}$
and $w_{n,0,n,k}$ go, and from $\id_{n,k}$ to $\id_{n,0}$ their
converses, $r^-_{n,k}$ and $w_{n,k,n,0}$. Between $\id_{n,i}$ and
$\id_{n,i}$ two atoms go, the identity atom $\id_{n,i}$ and the
diversity atom $w_{n,i,n,i}$. Between $\id_{n,i}$ and $\id_{m,j}$
when $n\ne m$ or $0\notin\{ i,j\}$ only one atom goes, namely
$w_{n,i,m,j}$.

The pair of the two diversity atoms $r_{n,k}$ and $w_{n,0,n,k}$
going from $\id_{n,0}$ to $\id_{n,k}$ is called a \emph{leaf}. These
leaves are connected at $\id_{n,0}$, and we call them a \emph{plant}
emanating from $\id_{n,0}$. So, intuitively, $\Z$ consists of
infinitely many plants, one plant with $n$ leaves for all natural
number $n$, see Figure~\ref{glove-fig}.

All the atoms of $\Z$ are ``big" in the sense that their products in
the complex algebra are as big as possible allowing
representability. At the end of this proof we show that the complex
algebra of $\Z$ is indeed completely representable, so it is a
relation algebra atom structure.

Let $\Z'$ be a nontrivial ultrapower of $\Z$. Then $\Z,\Z'$ are
elementarily equivalent by the fundamental theorem of ultraproducts.
We will show that their term algebras are not elementarily
equivalent, by exhibiting a formula that distinguishes them, as
follows.
The ultrapower $\Z'$ looks exactly like $\Z$, there are leaves
grouped into plants, except that in $\Z'$ there are plants with
infinitely many leaves, too. The key idea is that in the term
algebras, the finite and infinite plants can be distinguished by a
first-order logic formula. Namely, only finitely many leaves can be
``split" by any element of the term algebra. Thus, for an infinite
plant there is no element that splits all its leaves, while for a
finite plant clearly there is such an element. Using this, we can
express that in $\Z$ all the plants are finite while in $\Z'$ there
are infinite plants, too. This will be a first-order logic formula
distinguishing $\Z$ and $\Z'$. We note that in the complex algebra
of $\Z'$ there are elements splitting all the leaves of the infinite
plants, too, but these elements cannot be generated from singletons,
so they are not in the term algebra.

We now turn to elaborating the details of the above plan. Let $\Aa,
\Aa'$ denote the term algebras of $\Z, \Z'$, respectively. 
Let the operations of the term algebras be denoted as follows:
$\cdot,+,-,0,1$ denote Boolean intersection, addition,
complementation, zero and unit, and let $;,\conv, 1'$ denote
relation algebraic composition, conversion and identity element,
respectively. For $a,x,y\in A$ we say that \emph{$a$ splits $x,y$}
if $x,y$ are distinct identity atoms and there are atoms going from
$x$ to $y$ both in $a$ and in its complement. Formally, let
\begin{description}
\item{} $x\times y \de x;1\cdot 1;y $,
\item{} $\sigma(a,x,y)\deiff $\\
$ \mbox{$x,y$ are distinct identity atoms and $a\cdot (x\times y)\ne
0$, $-a\cdot (x\times y)\ne 0$}$.
\end{description}
Thus, if $a$ splits $x,y$, then $x,y$ is a splitable pair of
identity atoms. Next we prove by an easy induction that each element
in $A'$ splits only finitely many pairs of atoms. Let the induction
statement be denoted by $\P(a)$,
\begin{description}
\item{} $\P(a)\deiff \mbox{$a$ splits only finitely many pairs of
identity atoms}$.
\end{description}

Since $\Aa'$ is a term algebra, it is generated by its atoms. When
$a$ is an atom, it can split only one pair of atoms, namely it can
split only the pair $\dm{a},\rg{a}$ of its domain and range. Thus,
$\P(a)$ holds for all atoms $a$.

Assume that $\P(a), \P(b)$ hold. By definition, $a$ splits $x,y$ iff
$-a$ splits $x,y$, so $\P(-a)$ holds, too. It is easy to see that if
$a+b$ splits $x,y$ then at least one of $a$ and $b$ splits $x,y$. So
$a+b$ can split only those pairs of atoms that are split either by
$a$ or by $b$, hence $\P(a)$ and $\P(b)$ imply $\P(a+b)$.

$\P(1')$ holds, because $1'$ cannot split any pair of distinct atoms
(since $\dm{a}=\rg{a}$ for identity atoms).

If $a$ splits $x,y$, then $a\conv$ splits $y,x$, so if $\P(a)$
holds, then also $\P(a\conv)$ holds.

We set to show $\P(a;b)$.  First notice that the product $a;b$ of
two diversity atoms splits no pair of atoms in $\Aa$, and this fact
can be expressed by a first-order logic formula about $\Z$:
\begin{description}
\item{} $[a,b\notin I \land x,y\in I\land x\ne y\land P(a,b,c)\land
P(x,c,c)\land P(x,d,d)\land P(c,y,c)\land P(d,y,d)]\to P(a,b,d)$.
\end{description}
Therefore this same formula holds in $\Z'$, too, since $\Z'$ is an
ultrapower of $\Z$, and so the product $a;b$ of two diversity atoms
splits no pair of atoms in $\Aa'$, either. This implies
\begin{description}\label{div-eq}
\item{($\star$)}
 $\Aa'\models\forall a,b(a+b\le-1'\to \forall
x,y\lnot\sigma(a;b,x,y))$,
\end{description}
as follows. Let $x,y$ be a pair of distinct identity atoms, and let
$X,Y\subseteq S'$ be subsets of diversity atoms. Then $X;Y=\bigcup\{
a;b : a\in X, b\in Y\}$, in $\Aa'$. Now, assume that $x\times y\cap
X;Y\ne 0$. Then $x\times y\cap a;b\ne 0$ for some $a\in X, b\in Y$.
Since $a,b$ are diversity atoms, they do not split $x,y$, thus
$x\times y\subseteq a;b\subseteq X;Y$, which means that $X;Y$ does
not split $x,y$, either. This proves ($\star$).

Let now $a,b\in A'$ be arbitrary. Let $a'\de a\cdot 1'$, $a''\de
a\cdot -1'$, and the same for $b$, i.e., $b'\de b\cdot 1'$, $b''\de
b\cdot -1'$. Then $a;b=a';b + a;b' + a'';b''$. Hence, by the
previously proven case for addition, $a;b$ can split only those
pairs of atoms that are split by either one of $a';b, a;b'$ or
$a'';b''$. Assume that $x,y$ are distinct identity atoms. By the
above statement ($\star$), then $a'';b''$ does not split $x,y$. We
show that if $a';b$ splits $x,y$, then $b$ also splits them. Indeed,
let $r,s$ be atoms going from $x$ to $y$ such that $r\le a';b$ and
$s\le -(a';b)$. Then $r\le b$ since $a';b\le b$ by $a'\le 1'$. Also
we must have $x\le a'$. Assume $s\le b$. Then $s\le x;b\le a';b$,
contradicting our assumption $s\le-(a';b)$.
Similarly, if $a;b'$ splits $x,y$ then $a$ also splits them. Thus,
$a;b$ can only split those atoms that are split by $a$ or by $b$,
hence $\P(a;b)$ holds by $\P(a), \P(b)$. We have proved that $\P(a)$
holds for all $a\in A'$.

We are ready to exhibit the sentence distinguishing $\Aa$ and
$\Aa'$. This sentence states that for all atoms $x$ there is an
element $a$ that splits all pairs $x,y$ if they are splitable.
Formally, let
\begin{description}
\item{}
$\varphi\deiff \forall x\exists a\forall y(\exists
a\sigma(a,x,y)\to\sigma(a,x,y)).$
\end{description}
Now, $\Aa\models\varphi$ because for all $x$ there are only finitely
many $y$ such that $x,y$ is splitable. On the other hand, assume
that $\Z'$ is the ultrapower of $\Z$ by the ultrafilter $D$ on $J$,
and let
\begin{description}
\item{}
$x\de\langle \id_{n,0} : n\in J\rangle\slash D$.
\end{description}
Then $x$ is an identity atom in $\Aa'$ such that infinitely many
``leaves" emanate from $x$. Since all elements $a$ of $\Aa'$ can
split only finitely many of these leaves, there is no element in
$\Aa'$ that can split all of the splitable atoms $x,y$, thus
$\Aa'\not\models\varphi$. We have shown that $\Aa$ and $\Aa'$ are
not elementarily equivalent.

Finally, we show that $\Z$ and $\Z'$ are completely representable.
It is enough to show that $\Z$ is completely representable, since an
ultraproduct of completely representable atom structures is again
completely representable (see \cite[Exercise~14.1]{HH02}).

Let $U$ be a set and $R,T\subseteq U\times U$ be binary relations on
$U$. Relation composition $R\rp T$ of relations $R$ and $T$,
converse $R\inv$ of relation $R$ and the identity relation $\id_U$
on $U$ are defined as follows.
\begin{description}
\item{}
$R\rp T = \{ (u,v) : \exists w [(u,w)\in R\textit{ and }(w,v)\in T]
\}$,
\item{}
$R\inv = \{ (v,u) : (u,v)\in R \}$,
\item{}
$\id_U = \{ (u,u) : u\in U\}$.
\end{description}
The relation $U\times U\setminus\id_U$ is called the \emph{diversity
relation} on $U$.

To show complete representability of $\Z$, we will construct a set
$U$ and a function $\rep : S\to U\times U$ such that
\begin{description}
\item{(i)} $\langle \rep(a) : a\in S\rangle$ is a partition of
$U\times U$ to nonempty parts,
\item{(ii)} $a\in I$\quad iff\quad $\rep(a)\subseteq \id_U$\quad iff\quad $\rep(a)\cap\id_U\ne \emptyset$,
\item{(iii)} $C(a,b)$\quad iff\quad $\rep(b)=\rep(a)\inv$\quad iff\quad $\rep(b)\cap\rep(a)\inv\ne\emptyset$,
\item{(iv)} $P(a,b,c)$\ \ iff\ \ $\rep(c)\subseteq rep(a)\rp \rep(b)$\ \ iff\ \ $\rep(c)\cap(\rep(a)\rp\rep(b))\ne\emptyset$.
\end{description}
$U$ consists of countably many disjoint copies of $\omega$: let
$U_{n,i}= \omega\times \{(n,i)\}$, and
\[ U \de \bigcup\{ U_{n,i} : n\in\omega,  i\le n\} .\]
We define for $n,m\in\omega, i\le n, j\le m$
\begin{description}
\item{} $\rep(\id_{n,i})\de\{ (u,u) : u\in U_{n,i}\}$,\quad
$\rep(w_{n,i,n,i})\de U_{n,i}\times U_{n,i}\setminus\id_{n,i}$,
\item{} $\rep(w_{n,i,m,j})\de U_{n,i}\times U_{m,j}$\quad when $n\ne
m$, or $i\ne j$ and $0\notin\{ i, j\}$.
\end{description}
It remains to define $\rep(s)$ for the ``split" atoms $s$ (i.e., for
the atoms that form leaves). Assume $n,k\in\omega$ and $1\le k\le n$
(and $i,j\in\omega$ are arbitrary).
\begin{description}
\item{} $\rep(r_{n,k})\de$\\
$\{ \langle (i,n,0), (2^i3^j,n,k)\rangle : i,j\in\omega\}\cup
\{ \langle (i,n,0), (2^j3^i,n,k)\rangle : i,j\in\omega\}\cup$\\
$\{ \langle (2^i3^j,n,0), (i,n,k)\rangle : i,j\in\omega\}\cup
\{ \langle (2^j3^i,n,0), (i,n,k)\rangle : i,j\in\omega\}\cup$\\
$\{ \langle (2^i5^k,n,0), (i,n,k)\rangle : i,j\in\omega\}$.
\item{} $\rep(w_{n,0,n,k})\de U_{n,0}\times U_{n,k}\setminus
\rep(r_{n,k})$,
\item{} $\rep(r^-_{n,k})\de \rep(r_{n,k})\inv$, \quad $\rep(w_{n,k,n,0})\de \rep(w_{n,0,n,k})\inv$.
\end{description}
The properties of $R\de\rep(r_{n,k})$ that we will use are the
following. Let $T\de U_{n,0}\times U_{n,k}\setminus R$. Assume
$X,Y\in\{ R, T\}$ and $1\le\ell\le n$, $k\ne\ell$.
\begin{description}
\item{(r1)}  For each $u\in U_{n,0}$ there are at least two $v\in
U_{n,k}$ such that $(u,v)\in X$.
\item{(r2)} For each $v\in U_{n,k}$ there are at least two $u\in
U_{n,0}$ such that $(u,v)\in X$.
\item{(r3)} For all distinct $u,v\in U_{n,0}$ there is $w\in
U_{n,k}$ such that $(u,w)\in X$ and $(v,w)\in Y$.
\item{(r4)} For all distinct $u,v\in U_{n,k}$ there is $w\in
U_{n,0}$ such that $(w,u)\in X$ and $(w,v)\in Y$.
\item{(r5)} For all $u\in U_{n,k}$ and $v\in U_{n,\ell}$ there is $w\in
U_{n,0}$ such that $(w,u)\in X$ and $(w,v)\in Y$.
\end{description}
Now, (r1) follows from (r3) and (r2) follows from (r4). To check
$(r3)$, let $u,v\in U_{n,0}$ be distinct. Assume $u=(i,n,0)$ and
$v=(j,n,0)$. Then $i\ne j$ and $(u,w)\in R, (v,w)\in R$  for
$w=(2^i3^j,n,k)$ by the first line in the definition of $R$. Let
$w=(2^i3^q,n,k)$ where $q=j+1$. Then $(u,w)\in R$ by the first line
in the definition of $R$. To show $(v,w)\in T$ we have to show
$(v,w)\notin R$. This last statement is true because $j\notin\{
i,q\}$, so the pair $(v,w)$ is not included in $R$ by the first line
of its definition, and $j< 2^i3^q$ by $q=j+1$, so the pair $(v,w)$
is not included in $R$ by the second and third lines of its
definition. Similarly, $(u,w)\in T, (v,w)\in R$ for $w=(2^q3^j,n,k)$
where $q=i+1$. Finally, $(u,w)\in T, (v,w)\in T$ for
$w=(2^q3^t,n,k)$ where $q,t$ are both bigger than $i+j$. The proof
for (r4) is slightly more involved than the proof of (r3) because of
the last line in the definition of $R$. Let $u,v\in U_{n,k}$ be
distinct, say $u=(i,n,k)$ and $v=(j,n,k)$ with $i\ne j$. Then
$(w,u)\in R$ and $(w,v)\in R$ for $w=(2^i3^j,n,0)$ by the second
line in the definition of $R$. Let $w=(2^i3^q,n,k)$ where $q=j+1$.
Then $(w,u)\in R$ by the second line in the definition of $R$. Also,
the pair $(w,v)$ is not included in $R$ by the first line of its
definition since $2^i3^q>j$, it is not included in $R$ by the second
line since $j\notin\{ i,q\}$, and it is not included by the third
line since $5^k$ is not a divisor of $2^i3^q$. The rest of (r4) and
the case of (r5) when $i\ne j$ are similar. The third line in the
definition of $R$ is present for the case of (r5) when $i=j$: assume
$u=(i,n,k)$ and $v=(i,n,\ell)$, and let $w=(2^i5^k,n,0)$. Then
$(w,u)\in R$ by the third line in the definition of $R$ and
$(w,v)\notin R$ because $i<2^i5^k$ and $5^k$ is a divisor of
$2^i5^k$ but $5^{\ell}$ is not a divisor of $2^i5^k$.

We are ready to check (i)-(iv). To check (i)-(iii) is
straightforward. It is also straightforward to check (iv) when
$a,b,c$ are all diversity atoms that do not form leaves. It remains
to check the following when $k,\ell\ne 0, k\ne\ell$.
\begin{description}
\item{(iv.1)}
$\rep(w_{n,0,n,0})\rp\rep(a)\supseteq\rep(b)$\quad when $a,b\in\{
r_{n,k}, w_{n,0,n,k}\}$,
\item{(iv.2)}
$\rep(w_{n,k,n,k})\rp\rep(a)\supseteq\rep(b)$\quad when $a,b\in\{
r^-_{n,k}, w_{n,k,n,0}\}$,
\item{(iv.3)}
$\rep(w_{n,k,n,\ell})\rp\rep(a)\supseteq\rep(b)$\quad when $a\in\{
r^-_{n,\ell}, w_{n,\ell,n,0}\}$ and $b\in\{ r^-_{n,k},
w_{n,k,n,0}\}$,
\item{(iv.4)}
$\rep(a)\rp\rep(b)\supseteq\{ (u,v)\in U_{n,0}\times U_{n,0} : u\ne
v\}$\quad when $a\in\{ r_{n,k}, w_{n,0,n,k}\}$ and $b\in\{
r^-_{n,k}, w_{n,k,n,0}\}$,
\item{(iv.5)}
$\rep(a)\rp\rep(b)\supseteq\{ (u,v)\in U_{n,k}\times U_{n,k} : u\ne
v\}$\quad when $a\in\{ r^-_{n,k}, w_{n,k,n,0}\}$ and $b\in\{
r_{n,k}, w_{n,0,n,k}\}$,
\item{(iv.6)}
$\rep(a)\rp\rep(b)\supseteq U_{n,k}\times U_{n,\ell}$\quad when
$a\in\{ r^-_{n,k}, w_{n,k,n,0}\}$ and $b\in\{ r_{n,\ell},
w_{n,0,n,\ell}\}$.
\end{description}
Of the above, (iv.1) is true because of (r2) and because
$\rep(w_{n,0,n,0})$ is the diversity relation on $U_{n,0}$: assume
$(u,v)\in \rep(b)\subseteq U_{n,0}\times U_{n,k}$. There is
$(w,v)\in\rep(a), w\ne u$ because (r2). Then $w\in U_{n,0}$ because
the domain of $\rep(a)$ is $U_{n,0}$ and $(w,v)\in
\rep(w_{n,0,n,0})$ because $w\ne u$ and $\rep(w_{n,0,n,0})$ is the
diversity relation on $U_{n,0}$. Hence,
$(u,v)\in\{(u,w)\}\rp\{(w,v)\}\subseteq\rep(w_{n,0,n,0})\rp\rep(a)$.
The proofs of (iv.2) and (iv.3) are analogous. To show (iv.4), let
$u,v\in U_{n,0}$ be distinct. By (r3), there is $w\in U_{n,k}$ such
that $(u,w)\in\rep(a)$ and $(v,w)\in\rep(b)$, thus
$(u,v)\in\rep(a)\rp\rep(b)\inv$ (we used the already proven (iii)).
The rest is completely analogous, except that to prove (iv.5) we use
(r4), and to prove (iv.6) we use (r5).

By this, we have proved Theorem~\ref{main-thm}. \qed

\begin{rem}\label{splitting-rem}
{\rm We note that our atom structure is obtained from its
restriction to the atom structure where we omit the atoms $r_{n,k}$
and $r^-_{n,k}$ for all $n,k\in\omega, 1\le k\le n$. We get our $\Z$
by splitting $\id_{n,0}\times\id_{n,k}$ (and its converse) to get
the ``leaves". For the method of splitting in relation algebra see
\cite{AMaddNsplit}.}
\end{rem}

\begin{rem}\label{HHexample-rem}
{\rm We describe, briefly, another example of a pair of elementarily
equivalent relational algebra atom structures with non-elementarily
equivalent term algebras, due to Robin Hirsch and Ian Hodkinson.
(They devised this example after seeing ours \cite{HH17} and we
include the example with their permission.) Let $\HH=\langle
H,P,C,I\rangle$ be the atom structure where the atoms are the
elements of $H\de \omega\cup\{ 1',x,x^-\}$, there is one identity
atom, namely $I=\{ 1'\}$, all atoms are self-converse except for
$x,x^-$ which are each other's converses: $C=\{ (a,a) : a\in
H\setminus\{ x,x^-\}\}\cup\{ (x,x^-), (x^-,x)\}$, and there is one
kind of forbidden diversity triple, namely $(n,x,n+1)$ for
$n\in\omega$, i.e., $P=\{ (1',a,a) : a\in H\}\cup\{ (a,1',a) : a\in
H\}\cup\{ (a,b,1') : (a,b)\in C\}\cup\{ (a,b,c) : a,b,c\in
H\setminus\{ 1'\}, (a,b,c)\notin F\}$, where the set of the Peircean
transforms of the forbidden triples is $F\de \{ (n,x,n+1) :
n\in\omega\}\cup\{(n+1,x^-,n) : n\in\omega\}\cup\{(n+1,n,x^-) :
n\in\omega\}\cup\{(x^-,n,n+1) : n\in\omega\}\cup\{(n,n+1,x) :
n\in\omega\}\cup\{(x,n+1,n) : n\in\omega\}$. Let $\HH'$ be another
atom structure which is elementarily equivalent to $\HH$ but in
which there is a non-well-founded model $\omega^+$ instead of
$\omega$. In both term algebras the elements are finite and cofinite
sets of atoms. One can express that $n$ and $m$ are self-converse
diversity atoms and $m=n+1$, by using the forbidden triangles. Let
$\varphi$ be the formula that says that for all self-converse
non-identity atoms $n$ there is an element closed under predecessors
and containing $n$ but not $n+1$. This formula is true in the first
term algebra (the element is $\{ 0,1,2,\dots,n \}$) but false in the
second term algebra.}
\end{rem}

\bigskip\bigskip\bigskip

\noindent Alfr\'ed R\'enyi Institute of Mathematics, Hungarian
Academy of Sciences\\
Budapest, Re\'altanoda st.\ 13-15, H-1053 Hungary\\
andreka.hajnal@renyi.mta.hu, nemeti.istvan@renyi.mta.hu

\begin{thebibliography}{}

\bibitem{AMaddNsplit} Andr\'eka, H., Maddux, R.D., N\'emeti, I.,
Splitting in relation algebras. Proceedings of Amer. Math. Soc.
111,4 (1991) 1085--1093.

\bibitem{CAIII} Bezhanishvili, N., Varieties of two-dimensional cylindric
algebras. In: Cylindric-like algebras and algebraic logic., eds:
Andr\'eka, H., Ferenczi, M., N\'emeti, I., Bolyai Society
Mathematical Studies Vol 22, Springer Berlin Heidelberg New York
2013, pp.37--59.

\bibitem{G17} Givant, S., Advanced topics in relation algebras.
Springer International Publishing AG, 2017, xix + 605 pp.

\bibitem{HMT71} Henkin, L., Monk, J.D., Tarski, A., Cylindric
Algebras Part I, North-Holland, Amsterdam, 1971.

\bibitem{HH02} Hirsch, R., Hodkinson, I., Relation algebras by games. Studies in
Logic and the Foundations of Mathematics, vol. 147, Elsevier
Science, North-Holland Publishing Company, Amsterdam (2002)

\bibitem{HH17} Hirsch, R., Hodkinson, I., Letter to the authors.
October 23, 2017.

\bibitem{JT} J\'onsson, B., Tarski, A., Boolean algebras with operators.
Part I: Amer. J. Math. 73 (1951), 891-939. Part II. Amer. J. Math.
74 (1952), 127--162.

\bibitem{MTAMS} Maddux, R.D., Some varieties containing relation
algebras. Trans. Amer. Math. Soc. 272,2 (1982), 501--526.

\bibitem{Mbook} Maddux, R.D., Relation algebras. Studies in Logic and the
Foundations of Mathematics, vol. 150. Elsevier Science,
North-Holland Publishing Company, Amsterdam (2006)

\bibitem{SDis} Simon, A., Non-representable algebras of relations. PhD Dissertation,
Hungarian Academy of Sciences, Budapest, 1997. iii+86pp.

\bibitem{Ven} Venema, Y., Atom structures. In: Kracht, M., de Rijke, M.,
Wansing, H., Zakharyaschev, M., eds., Advances in Modal Logic 96.
pp.291-305. CSLI Publications, Stanford, 1997.
\end{thebibliography}
\end{document}